\documentclass[12pt]{article}
\usepackage{latexsym}
\usepackage{amsfonts}
\usepackage{amsmath}
\usepackage{amssymb}
\usepackage{graphicx}%
\def\demo{\noindent{\bf Proof. }}
\def\QED{\hfill$\Box$}

\newtheorem{Theorem}{Theorem}[section]
\newtheorem{Lemma}[Theorem]{Lemma}
\newtheorem{Corollary}[Theorem]{Corollary}
\newtheorem{Proposition}[Theorem]{Proposition}
\newtheorem{Remark}[Theorem]{Remark}

\newtheorem{Conjecture}[Theorem]{Conjecture}
\newtheorem{Definition}[Theorem]{Definition}

\begin{document}
\topmargin3mm
\hoffset=-1cm
\voffset=-1.5cm
\

\medskip

\begin{center}
{\large\bf On 2-partitionable clutters and the MFMC property}
\vspace{6mm}\\
\footnotetext{2000 {\it Mathematics Subject
Classification}. Primary 05C75; Secondary 05C85, 05C20,
13H10.}
\footnotetext{{\it Key words and phrases\/}.
max-flow min-cut, mengerian, clutters, hypergraphs, normality, Rees algebras}

\medskip

Alejandro Flores-M\'endez, \footnote{LIDETEA Universidad La Salle}
Isidoro Gitler \footnote{This work was partially supported
by CONACyT grants 49251-F, 49835-F, and SNI.}
and Enrique Reyes
\footnote{Partially supported by COFAA-IPN.}
\\
{\small Departamento de Matem\'aticas}\vspace{-1mm}\\
{\small Centro de Investigaci\'on y de Estudios Avanzados del
IPN}\vspace{-1mm}\\
{\small Apartado Postal 14--740}\vspace{-1mm}\\
{\small 07000 M\'exico City, D.F.}\vspace{-1mm}\\
{\small e-mail: {\{{\tt aflores,igitler,ereyes}\}{\tt @math.cinvestav.mx}}}\vspace{4mm}
\end{center}
\date{}

\begin{abstract}
\noindent

We introduce 2-partitionable clutters as the simplest case of the class of $k$-partitionable 
clutters and study some of their combinatorial properties. In particular, we study properties of the rank of
the incidence matrix of these clutters and properties of their minors.

A well known conjecture of Conforti and Cornu\'ejols ~\cite{ConfortiCornuejols,cornu-book} states: That 
all the clutters with the packing property have the max-flow min-cut property, i.e. are mengerian.
Among the general classes of clutters known to verify the conjecture
are: balanced clutters (Fulkerson, Hoffman and
Oppenheim~\cite{FulkersonHoffmanOppenheim}), binary clutters
(Seymour~\cite{Seymour}) and dyadic clutters (Cornu\'ejols, Guenin
and Margot~\cite{CornuejolsGueninMargot}). We find a new infinite
family of 2-partitionable clutters, that verifies the conjecture.

On the other hand we are interested in studying the normality of the Rees algebra associated to a
clutter and possible relations with the Conforti and Cornu\'ejols conjecture.
In fact this conjecture is equivalent to an algebraic statement about the normality
of the Rees algebra~\cite{rocky}.

\end{abstract}

\medskip

\section{Introduction}
We briefly describe the main results in this paper.
Theorem ~\ref{menger} characterizes when an ideal clutter is mengerian in terms of the existence 
of an edge $e$ of $H$ for which $\tau^{\omega}\left( H \right) = \tau^{\omega-e} \left( H \right) + 1$,
for all $\tau^{\omega}\left( H \right) > 0$. Hence the Conforti-Cornu\'ejols conjecture reduces to 
proving that for every hypergraph $H$ with the packing property if $\tau^{\omega}\left( H \right)>0$, 
for some $\omega$, then there exists $e\in E(H)$, for which $\tau^{\omega}\left( H \right) = \tau^{\omega-e} \left( H \right) + 1$.

After introducing and showing that the family of hypergraphs $Q_{pq}^{F}$ has the packing property,
we use Theorem ~\ref{menger} to prove that it is also mengerian. We give an explicit algorithm
in pseudo code to obtain inductively the edge $e$ required in Theorem ~\ref{menger}.
This algorithm can be generalized to prove that other ideal hypergraphs are mengerian.

We introduce 2-partitionable clutters. We then prove Propositions ~\ref{jan17-04} and ~\ref{2-part}
that give information on the rank of the incidence matrix $A$ of the clutter $H$ and the minors of $H$ (when $H$ is 2-partitionable). 
We propose the Conforti and Cornuejos conjecture
for 2-partitionable hypergraphs since we believe that in the class of $k$-partitionable hypergraphs a
counterexample can in principle be possible.

In proposition ~\ref{Deltar-Q6} we prove the following conjecture 
for the case of the hypergraphs $Q_{pq}^{F}$, thus giving support to this conjecture.

Conjecture~\ref{git-val-vi}: If $\tau({\cal C}')=\nu({\cal C}')$ for
all minors $\cal C'$ of $\cal C$ and $x^{v_1},\dots,x^{v_q}$
have degree $d\geq 2$, then the group
$
\mathbb{Z}^{n+1}/((v_1,1),\ldots,(v_q,1))
$
is free, or equivalently $\Delta_r(B)=1$ where $r={\rm rank}(B)$.

\section{Preliminaries}
We now give several definitions that help to clarify the algebraic translations
of the combinatorial optimization problems studied in this paper.

A {\it hypergraph} $H$ is defined by a pair $\left(  V,E\right)$ where $V$
represents a finite set called the vertices of $H$, and $E$
represents a collection of subsets from $V$ called the edges of $H$.
In some cases, the notation $V\left(  H\right)  , E\left(
H\right)  $ will be used to refer to the vertices and the edges of $H$, respectively.

A {\it clutter\/} $\cal C$ is a particular type of hypergraph, with
the property that $S_1 \not\subseteq S_2$ for all distinct $S_1,S_2
\in E(\cal C)$.
%We view an hypergraph $H$ as a clutter if for every pair of edges $e$, $f$ such that %$e\subseteq f$ this implies that $e=f$. A basic example of clutter is a graph.

For every hypergraph $H=\left(
V,E\right)  $ there is an {\it associated clutter} $H^{\min}$, defined by:%
\[
H^{\min}=\left\{  e\in E:e\nsupseteq f\in E\right\}
\]

The contraction $H / i$ and deletion $H \backslash i$ are
hypergraphs with vertex set $V( H )\setminus \{i\}$ where: $E(H /
i)=\{ S \setminus \{ i \} :S \in E(H) \}$ (for clutters we take the
set of inclusionwise minimal members of this set) and $E(H
\backslash i) = \{ S \in E(H) : i \notin S \}$. Contractions and
deletions of distinct vertices can be performed sequentially and the
result does not depend on the order. An hypergraph obtained from $H$
by a sequence of deletions $I_d$ and contractions $I_c$  $(I_d \cap
I_c = \emptyset)$ is called a minor of $H$ and is denoted by $H
\backslash I_d  / I_c$. If $I_d \not= \emptyset$ or $I_c \not=
\emptyset$, the minor is proper. For general properties of
hypergraphs, clutters and their blockers we refer the reader
to~\cite{Schrijver}.

Let $R=K[x_1,\ldots,x_n]$ be a polynomial ring
over a field $K$ and let $I$ be an ideal
of $R$ of height $g\geq 2$, minimally generated by a finite set
$F=\{x^{v_1},\ldots,x^{v_q}\}$
of square-free monomials of degree at least two. Where a monomial $f$ in $R$ is called
{\it square-free\/} if
$f=x_{i_1}\ldots x_{i_r}$
for some $1\leq i_1<\cdots<i_r\leq n$. For technical reasons we
shall assume that each variable $x_i$ occurs in at least one
monomial of $F$.

There is a natural one to one
correspondence between the family of
square-free monomial ideals and the family of clutters:

We associate to the
ideal $I$ a {\it clutter\/} $\cal C$ by taking the set
of indeterminates $V=\{x_1,\ldots,x_n\}$ as vertex set and
$E=\{S_1,\ldots,S_q\}$ as edge set, where
$$S_k=\{x_i | \, \langle e_i,v_k\rangle=1, e_i \,\mbox{denotes the i'th unit vector} \}={\rm supp}(x^{v_k}),$$
and the {\it support\/} of a monomial $x^a=x_1^{a_1}\cdots
x_n^{a_n}$ in $R$ is given by ${\rm supp}(x^a)= \{x_i\, |\,
a_i>0\}$. The ideal $I$ is called the {\it edge ideal\/} of $\cal
C$. To stress the relationship between $I$ and $\cal C$ we will use
the notation $I=I({\cal C})$. The $\{0,1\}$-vector $v_k$ is called
the {\it characteristic vector\/}\index{characteristic vector} of
$S_k$. We denote by  ${\mathbf 1}$ the vector whose entries are all
ones.

Let $A$ be the {\it incidence matrix\/} whose column vectors are
$v_1,\ldots,v_q$.
%In order to link the algebraic properties of these algebras with combinatorial
%optimization
The {\it set covering polyhedron\/} associated with $A$ is defined as:

$$
Q(A)=\{x\in\mathbb{R}^n |\, x\geq 0;\, xA\geq \mathbf{1}\}.
$$

We say that a hypergraph (clutter) is \emph{ideal} if the polyhedron $Q(A)$ is integral.

\medskip

A set $C\subseteq V$ is a {\it vertex cover or transversal \/} of
the hypergraph $H$ if every edge of $H$ contains at least one vertex
in $C$. We call $C$ a {\it minimal vertex cover or minimal
transversal\/} if $C$ is minimal with respect to this property. A
set of edges of the hypergraph $H$ is called {\it independent or a
matching\/} if no two of them have a vertex in common.

%Let $H=\left(  V,E\right)  $ and
%$t\subseteq V$, recall that $t$ is a
%transversal of $H$ if $t\cap e\neq\emptyset$ for every $e\in E$.
Let $B\left(  H\right)  $ be the collection of all the transversals
of $H$, then
the \emph{blocker} of $H$, $b\left(  H\right)  $, is defined by:%
\[
b\left(  H\right)  :=B\left(  H\right)  ^{\min}%
\]

It is well known that the following dual relationship holds for
the blocker.

\begin{Theorem}\label{bloker}
For every hypergraph $H$, $b\left(  b\left(  H\right)  \right)  =H^{\min}$. In
particular, if $H$ is a clutter, then $H=b\left(  b\left(  H\right)  \right)
$.
\end{Theorem}

Observe that for minors we have: $b(H \backslash i) = (b(H) /
i)^{\min}$ and $b(H / i) = b(H)\backslash i$.

\medskip

Let us denote by $f^c$ the complement of $f$. Then as an immediate consequence of the
former theorem we have that.

\begin{Corollary}
\label{f_Vmenosf}For every hypergraph $H=\left(  V,E\right)  $, and
$f\subseteq V$, either $f\supseteq e\in E$, or $f^c \supseteq t\in
b\left(  H\right)  $, but not both.
\end{Corollary}

\begin{Proposition}\label{1cover-1}{\rm \cite{rocky}} The following are equivalent\/{\rm :}
\begin{description}
\item{\rm(a)} $\mathfrak{p}=(x_{1},\ldots,x_{r})$ is a minimal prime of
$I=I({\cal C})$.
\item{\rm(b)} $C=\{x_{1},\ldots,x_{r}\} \in b({\cal C})$.
\item{\rm(c)} $\alpha=e_{1}+\cdots+e_{r}$ is a vertex of $Q(A)$.
\end{description}
\end{Proposition}

\subsection{The Conforti and Cornu\'ejols conjecture}

\begin{Definition}\rm The clutter $\cal C$ satisfies the
{\it max-flow min-cut\/}
(MFMC)
property if both sides
of the LP-duality equation
\begin{equation}\label{jun6-2-03-1}
\tau^{\omega}={\rm min}\{\langle \omega,x\rangle |\, x\geq 0; xA\geq \mathbf{1}\}=
{\rm max}\{\langle y,\mathbf{1}\rangle |\, y\geq 0; Ay\leq\omega\}=\nu^{\omega}
\end{equation}
have integral optimum solutions $x$ and $y$
for each non-negative integral vector $\omega$. Clutters that satisfy the MFMC property are called \emph{mengerian}.
\end{Definition}

 \medskip

 Recall that a monomial
subring $K[F]\subset R$ is {\it normal\/} if
$K[F]=\overline{K[F]}$, where the
integral closure $\overline{K[F]}$ is given by:
\begin{equation}\label{jan5-02}
\overline{K[F]}=K[\{x^a\vert\, a\in
{\mathbb Z}{\cal A}\cap {\mathbb R}_+{\cal A}\}],
\end{equation}
where $\mathbb{Z}{\cal A}$ is the subgroup spanned by $\cal A$
and $\mathbb{R}_+{\cal A}$ is the {\it polyhedral cone\/}
\[
{\mathbb R}_+{\cal A}= \left.\left\{\sum_{i=1}^qa_i v_i\right\vert\,
a_i\in {\mathbb R}_+\, \mbox{ for all }i\right\}
\]
generated by ${\cal A}=\{v_1,\ldots,v_q\}$. Here $\mathbb{R}_+$
denotes the set of non negative real numbers.

\begin{Theorem}\label{noclu1}{\rm \cite{normali,clutters,HuSV}} The following are equivalent
\begin{description}
 %\item{\rm(i)\ \ } ${\rm gr}_I(R)$ is reduced.
 %\vspace{-1mm}
 \item{\rm (i)} $\cal C$ has the max-flow min-cut property.
 \item{\rm (ii)\ } $R[It]$ is normal and $Q(A)$ is an integral
 polyhedron, where $R[It]$ is the
subring of $R[t]$ generated by $x_1,\ldots ,x_n,
x^{v_1}t,\ldots,x^{v_q}t$, over the field $K$.\vspace{-1mm}
 \item{\rm (iii)} $x\geq 0;\, xA\geq \mathbf{1}$ is totally dual integral ({\rm a TDI system}).

\end{description}
\end{Theorem}

\begin{Proposition}\label{propsey}{\rm \cite{Seymour}} If a clutter
$\cal C$ has the MFMC property, then so do all its minors.
\end{Proposition}

\medskip

Let us denote by $\tau({\cal C})$ the minimum size of a vertex cover in $\cal C$ and
by $\nu({\cal C})$ the maximum size of a matching in ${\cal C}$, then:

\begin{eqnarray*}
 \lefteqn{\tau({\cal C})\geq {\rm min}\{\langle \mathbf{1},x\rangle |\, x\geq 0; xA\geq
 \mathbf{1}\}}\\
 &\ \ \ \ \ \ \ \ \ \ &
 ={\rm max}\{\langle y,\mathbf{1}\rangle |\, y\geq 0; Ay\leq\mathbf{1}\}
 \geq \nu({\cal C}).
 \end{eqnarray*}

%Note that $\tau({\cal C})=\nu({\cal C})$ if and only if both sides of
%the equality have integral optimum solutions. These two numbers can
%be interpreted
%in terms of invariants of $I$.
%The height
%of the ideal $I$, is denoted by ${\rm ht}(I)$. We have for the {\it vertex covering number} that: $\tau({\cal C})={\rm ht}(I)$.
%
%
%On the
%other hand the {\it edge independence number\/} can be expressed as:
%$\nu({\cal C})=\max\{r|\, \exists\, \mbox{ a regular sequence of monomials }
%x^{\alpha_1},\ldots,x^{\alpha_r}\in I\}.$
%The number in the right hand side is denoted by
%${\rm mgrade}(I)$, it is called
%the {\it monomial grade\/} of the ideal.

\begin{Definition} A clutter $\cal C$ has the K\"onig property or packs if
$\tau({\cal C})=\nu({\cal C})$.
\end{Definition}

The simplest example of a clutter with the K\"onig property
is a bipartite graph.

\medskip

\begin{Definition}\rm A clutter $\cal C$ satisfies the {\it packing property\/}
(PP) if all its minors satisfy the K\"onig property,
that is, $\tau({\cal C}')=\nu({\cal C}')$
for every minor ${\cal C}'$ of $\cal C$.
\end{Definition}

%\begin{Proposition} Let $x^{z_1},\ldots,x^{z_r}$ be a set of
%generators of $I$ and let $B$ be the matrix whose columns are
%$z_1,\ldots,z_r$.
%Then both sides of the equation
%$$
%{\rm min}\{\langle \mathbf{1},x\rangle |\, x\geq 0; xA\geq
%\mathbf{1}\}={\rm max}\{\langle y,\mathbf{1}\rangle |\, y\geq 0; Ay\leq\mathbf{1}\}
%$$
%have integral
%optimum solutions if and only if both sides of the equation
%$$
%{\rm min}\{\langle \mathbf{1},x\rangle |\, x\geq 0; xB\geq
%\mathbf{1}\}={\rm max}\{\langle y,\mathbf{1}\rangle |\, y\geq 0; By\leq\mathbf{1}\}
%$$
%have integral optimum solutions. If any of the two systems have
%integral optimum solutions $x,y$ their optimal values are equal.
%\end{Proposition}
%
%\demo
%The idea is to show that the integrality conditions
%depend only on the monomial ideal $I$. From the inequalities:
%\begin{eqnarray*}
%\lefteqn{{\rm ht}(I)\geq {\rm min}\{\langle \mathbf{1},x\rangle
%|\, x\geq 0; xB\geq
%\mathbf{1}\}}\\
%&\ \ \ \ \ \ \ \ \ \ &
%={\rm max}\{\langle y,\mathbf{1}\rangle |\, y\geq 0; By\leq\mathbf{1}\}
%\geq {\rm mgrade}(I)
%\end{eqnarray*}
%we rapidly obtain that ${\rm ht}(I)={\rm mgrade}(I)$ if and only
%if both sides of the min-max equation have integral optimum solutions
%$x,y$. Therefore the assertion follows by observing that the columns of $A$
%correspond to a minimal set of generators of the ideal $I$.
%\QED

\begin{Theorem}\label{lehman}{\rm \cite{Lehman}} If\, $\cal C$ has the packing property, then $Q(A)$ is
integral.
\end{Theorem}

It is well known that:

\begin{Proposition}\label{jan20-05} If $\cal C$ has the max-flow min-cut property,
then $\cal C$ has the packing property.
\end{Proposition}

\begin{Conjecture}\label{conforti-cornuejols1}{\rm(Conforti-Cornu\'ejols~\cite{ConfortiCornuejols,cornu-book})}
\rm  If the clutter $\cal C$ has the
packing property, then $\cal C$ has the max-flow min-cut property.
\end{Conjecture}

Next we state an algebraic
version of Conjecture~\ref{conforti-cornuejols1}.

\begin{Conjecture}\label{con-cor-vila}{\rm \cite{rocky}} \rm If $\tau({\cal C}')=\nu({\cal C}')$ for
all minors $\cal C'$ of $\cal C$, then $R[It]$ is normal.
\end{Conjecture}

\noindent {\it Notation} For an integral matrix $B\neq(0)$, the
greatest common divisor of all the nonzero $r\times r$
subdeterminants of $B$ will be
denoted by $\Delta_r(B)$.

\begin{Theorem}\label{MM}{\rm \cite{rocky}} If $x^{v_1},\ldots,x^{v_q}$
are monomials of degree $d\geq 2$, i.e., all the edges of the clutter have $d$
vertices, such that $\cal C$ satisfies MFMC and the matrix
$$
B=\left(\hspace{-1mm}
\begin{array}{ccc}
v_1&\cdots&v_q\\
1&\cdots &1
\end{array}\hspace{-1mm}
\right)
$$
has rank $r$, then $\Delta_r(B)=1$.
\end{Theorem}

By using the previous result we obtain that a positive answer to
Conjecture~\ref{con-cor-vila}  implies the following:

\begin{Conjecture}\label{git-val-vi}{\rm \cite{rocky}} \rm\ If $\tau({\cal C}')=\nu({\cal C}')$ for
all minors $\cal C'$ of $\cal C$ and $x^{v_1},\dots,x^{v_q}$
have degree $d\geq 2$, then the group
$$
\mathbb{Z}^{n+1}/((v_1,1),\ldots,(v_q,1))
$$
is free, or equivalently $\Delta_r(B)=1$ where $r={\rm rank}(B)$
\end{Conjecture}

\section{On 2-partitionable clutters}

%Take $R=K[x_1,\ldots,x_n]$
%be a polynomial ring
%over a field $K$
%and let
%$$V=X_1\cup X_2\cup\cdots\cup X_d$$
%be a partition of $V=\{x_1,\ldots,x_n\}$ into $d$ subsets of size
%two. We set $X_i=\{x_{2i-1},x_{2i}\}$ for $i=1,\ldots,d$ and $d\geq 2$. Note
%$n=2d$.
%
%For ${\cal C}$ 2-partitionable and $I=I({\cal C})$ minimally generated by
%$F=\{x^{v_1},\ldots,x^{v_q}\}$ we have that
%\begin{equation}\label{jan19-04}
%|{\rm supp}(x^{v_i})\cap X_k|=1\ \ \ \ \forall\ i,k.
%\end{equation}

Let $H=(V,E)$ be a hypergraph with $|V|=2d$ and $E= \{ S_1, \ldots , S_k \} $. Let $$V=X_1\cup X_2\cup\cdots\cup X_d$$
be a partition of $V=\{x_1,\ldots,x_n\}$ into $d$ subsets of size
two. We set $X_i=\{x_{2i-1},x_{2i}\}$ for $i=1,\ldots,d$ and $d\geq 2$.
If
\begin{equation}\label{jan19-04}
|S_i\cap X_j|=1\ \ \ \ \forall\ i,j.
\end{equation}

\noindent we say that $H$ is \emph{$2$-partitionable}. This definition could be generalized to 
\emph{$k$-partitionable} hypergraphs, where $|X_j|=k\geq 2$.

Note that for $I=I({\cal C})$ minimally generated by
$F=\{x^{v_1},\ldots,x^{v_q}\}$ we have that equation (\ref{jan19-04}) becomes:
\begin{equation}\label{jan19-06}
|{\rm supp}(x^{v_i})\cap X_j|=1\ \ \ \ \forall\ i,j.
\end{equation}

%As before, we assume that each variable $x_i$ occurs in at least one
%monomial of $F$. In the sequel $A$ will denote the $n\times q$
%incidence matrix with column vectors $v_1,\ldots,v_q$, and $\cal C$
%will denote the clutter associated to $I$.

%In this case the clutter
%${\cal C}$ is called \emph{2-partitionable}.

\medskip

Observation: In our situation, by the pigeon hole
principle, any minimal vertex cover $C$ of the clutter $\cal C$
satisfies
$2\leq|C|\leq d$. Notice
that for each odd integer $k$ the sum of rows $k$ and $k+1$ of the matrix
$A$ is equal to ${\mathbf 1}=(1,\ldots,1)$. Thus the rank of $A$ is
bounded by $d+1$.

The next result
shows that $A$ has ``maximal rank'' if $\cal C$ has a cover of
maximum possible size.

\begin{Proposition}\label{jan17-04} Let $\cal C$ be a 2-partitionable
clutter. If there exists a minimal vertex cover $C$ such that $|C|=d\geq 3$ and $\cal C$ satisfies the K\"onig property, then ${\rm rank}(A)=d+1$.
\end{Proposition}

\demo First notice that $C$ contains exactly one element of each
$X_j$ because $X_j\not\subset C$. Thus we may assume
$$
C=\{x_1,x_3,\ldots,x_{2d-1}\}.
$$
Consider the monomial $x^\alpha=x_2x_4\cdots x_{2d}$
and notice that $x_kx^\alpha\in I$ for each $x_k\in C$ because the
monomial $x_kx^\alpha$ is clearly in every minimal prime of $I$.
Writing $x_k=x_{2i-1}$ with $1\leq i\leq d$ we conclude that the
monomial
$$
x^{\alpha_i}=x_2x_4\cdots x_{2(i-1)}x_{2i-1}x_{2(i+1)}\cdots x_{2d}
$$
is a minimal generator of $I$. Thus we may assume
$x^{\alpha_i}=x^{v_i}$ for $i=1,\ldots,d$. The vector $\mathbf 1$
belongs to the linear space generated by $v_1,\ldots,v_q$ because
$\cal C$ has the K\"onig property. It follows readily that the
matrix with rows $v_1,\ldots,v_d,\mathbf{1}$ has rank $d+1$. \QED

\medskip

\begin{Remark}\label{RR}
If ${\cal C}$ is $2$-partitionable then $rank(A)=rank(B)=d+1$, where
$$
B=\left(\hspace{-1mm}
\begin{array}{ccc}
v_1&\cdots&v_q\\
1&\cdots &1
\end{array}\hspace{-1mm}
\right)
$$
\end{Remark}

\begin{Proposition}\label{2-part}
Let $H$ be $2-$partitionable hypergraph with $W,Z\subseteq V\left(
H\right) $ such that $W\cap Z=\emptyset$ and $Z\cap X_{i}=\emptyset$
for some $X_{i}$. If $H^{\prime}\backslash W/Z$ is a proper minor of
$H$, then $H^{\prime}$ has the K\"onig property or $W\cap
X_{i}=\emptyset$.

\end{Proposition}

\demo If $\emptyset\in E\left(  H^{\prime}\right)  $ or $E\left(
H^{\prime}\right) =\emptyset$, then $H^{\prime}$ has the K\"onig
property. We assume that $\emptyset\notin E\left( H^{\prime}\right)
$ and $E\left( H^{\prime}\right)  \neq\emptyset$. If $\tau\left(
H^{\prime}\right) =1$ then $H^{\prime}$ has the K\"onig property.
Let us consider then that $\tau\left( H^{\prime}\right) \geq2$, and
$X_{i}=X_{1}=\left\{
x_{1},x_{2}\right\}  $. As $Z\cap X_{1}=\emptyset$ then $\left\{  x_{1}%
,x_{2}\right\}  \in b\left(  H/Z\right)  =b\left(  H\right)
\backslash Z$. Hence, $\left\vert \left\{  x_{1},x_{2}\right\}
\backslash W\right\vert \geq\tau\left(  H^{\prime}\right)  \geq2$
and therefore $W\cap X_{1}=\emptyset$. \QED

\begin{Conjecture}\label{2part}
Let ${\cal C}$ be a 2-partitionable clutter. Then ${\cal C}$ has the
packing property if and only if it is mengerian.
\end{Conjecture}

\section{The $Q_6$ property class of hypergraphs}

A hypergraph $H$ is minimally non packing (MNP) if it does not pack,
but every minor of it does. Cornu\'ejols, Guenin and
Margot~\cite{CornuejolsGueninMargot}, give an infinite class of
ideal MNP clutters, which they call the $Q_{6}$ property class
(before their work, only two MNP clutters were known).

\medskip
A clutter has the $Q_6$ property, when $V({\cal C})$ can be partitioned into
nonempty sets $I_1, \ldots , I_6$ such that there are edges $S_1, \ldots , S_4$ in
${\cal C}$ of the form:

$$S_1=I_1 \cup I_3 \cup I_5 , \quad \quad S_2=I_1 \cup I_4 \cup I_6 , \quad \quad S_3=I_2 \cup I_4 \cup I_5 , \quad \quad S_4=I_2 \cup I_3 \cup I_6 .$$

\medskip

The Cornu\'ejols, Guenin and Margot MNP family of $Q_6$-property
clutters is described as follows. Given $p,q\in\mathbb{N}$, we
construct the incidence matrix of the clutter $Q_{pq}$ by
partitioning the set $V\left(  Q_{pq}\right)$ in 6 blocks which we
will call $P,P^{\ast},Q,Q^{\ast},r,r^{\ast}$, with elements
$P=\left\{ p_{1},\ldots,p_{p}\right\},$ $P^{\ast}=\left\{
p_{1}^{\ast},\ldots,p_{p}^{\ast}\right\},$ $Q=\left\{
q_{1},\ldots,q_{q}\right\},$
$Q^{\ast}=\left\{q_{1}^{\ast},\ldots,q_{q}^{\ast}\right\}$.
Furthermore, denote by
$M_{m \times n}\left(  \mathbb{B}%
\right)  $ the set of $0,1$ matrices and let $H_{p}\in
M_{\left(  \left(  2^{p}-1\right)  \times p\right)  }\left(  \mathbb{B}%
\right)  $ be a matrix whose rows represent the characteristic vectors of the non empty
subsets of a set with $p$ elements. Let $H_{p}^{\ast}$ be its complement,
i.e. $H_{p}+H_{p}^{\ast}=J$, where $J$ denotes the matrix whose entries are
all one. Then the transpose $A^t$ of the incidence matrix $A$ of the clutter $Q_{pq}$ is given by:%
\[%
\begin{tabular*}{3.2in}
[r]{cp{0.5in}p{0.5in}p{0.5in}p{0.55in}p{0.25in}p{0.25in}}
& $\, p_{1}\ldots p_{p}$ & $p_{1}^{\ast}\ldots p_{p}^{\ast}$ & $q_{1}\ldots
q_{q}$ & $q_{1}^{\ast}\ldots q_{q}^{\ast}$ & $r$ & $r^{\ast}$\\
\multicolumn{1}{r}{$A^t\left(  Q_{pq}\right)  =$} & \multicolumn{6}{l}{$\left[
\begin{tabular*}{3.2in}
[c]{p{0.5in}p{0.5in}p{0.5in}p{0.35in}p{0.25in}p{0.25in}}%
$H_{p}$ & $H_{p}^{\ast}$ & $J$ & $\mathbf{0}$ & $\mathbf{1}$ & $\mathbf{0}$\\
$H_{p}^{\ast}$ & $H_{p}$ & $\mathbf{0}$ & $J$ & $\mathbf{1}$ & $\mathbf{0}$\\
$J$ & $\mathbf{0}$ & $H_{q}^{\ast}$ & $H_{q}$ & $\mathbf{0}$ & $\mathbf{1}$\\
$\mathbf{0}$ & $J$ & $H_{q}$ & $H_{q}^{\ast}$ & $\mathbf{0}$ & $\mathbf{1}$%
\end{tabular*}
\right]  $}%
\end{tabular*}
\]

The hypergraph $Q_{6}$ giving name to the class, corresponds to $Q_{1,1}$.

\medskip

As an example of an hypergraph in this class, we show the incidence matrix of $Q_{2,1}$:%
\[
A^t\left(  Q_{2,1}\right)  =%
\begin{bmatrix}
1 & 1 & 0 & 0 & 1 & 0 & 1 & 0\\
1 & 0 & 0 & 1 & 1 & 0 & 1 & 0\\
0 & 1 & 1 & 0 & 1 & 0 & 1 & 0\\
0 & 0 & 1 & 1 & 0 & 1 & 1 & 0\\
0 & 1 & 1 & 0 & 0 & 1 & 1 & 0\\
1 & 0 & 0 & 1 & 0 & 1 & 1 & 0\\
1 & 1 & 0 & 0 & 0 & 1 & 0 & 1\\
0 & 0 & 1 & 1 & 1 & 0 & 0 & 1
\end{bmatrix}
\]

\subsection{The family $Q_{pq}^{F}$ of $2$-partitionable hypergraphs}

For homogeneity reasons, from here on we will assume that $p,q>1$. From the
construction of $Q_{pq}$, it follows that every pair of vertices of the form $vv^{\ast}$
with $v\in PQr$ is contained in $b\left(  Q_{pq}\right)  $ ($XYz$ will be used
as a shorthand for the union of sets $X\cup Y\cup\left\{  z\right\}  $). Moreover, it is not
hard to show that the elements of $b\left(  Q_{pq}\right)  $ correspond to one
of the following types:
\begin{subequations}
\label{eq_bQpq}%
\begin{gather}
vv^{\ast}\text{, where }v\in PQr\\
p_{i}p_{j}^{\ast}q_{k}q_{l}^{\ast}\text{, where }i\neq j,k\neq l\\
p_{i}p_{j}^{\ast}r\text{, where }i\neq j\\
q_{i}q_{j}^{\ast}r^{\ast}\text{, where }i\neq j\\
Pq_{i}^{\ast}r^{\ast}\\
P^{\ast}q_{i}r^{\ast}\\
p_{i}Qr\\
p_{i}^{\ast}Q^{\ast}r
\end{gather}
\end{subequations}

For simplicity, we will denote by $V_{pq},E_{pq}$ the set of
vertices $V\left(  Q_{pq}\right)$ and the set of edges $E\left(
Q_{pq}\right)$, respectively. If $A\in \{P,Q\}$ and $A_1\subseteq A$
we define $A_{1}^{c}=A \backslash A_1$. We will also use the
following definitions:

\begin{align*}
F_{pq}  &  :=\left\{  PQr,P^{\ast}Q^{\ast}r,PQ^{\ast}r^{\ast},P^{\ast}%
Qr^{\ast}\right\} \\
Q_{pq}^{F}  &  :=\left(  V_{pq},E_{pq}\cup F\right),
\text{where} \, F\subseteq E_{pq}^{\ast} \\
F^{\ast}  &  :=\left\{  f^{\ast}:f\in F\right\} \, \text{ note that } \,
f=\left(  f^{\ast}\right)  ^{\ast} \\
\end{align*}

\medskip

It is easy to verify that by construction the hypergraphs $Q_{pq}^{F}$ are $2$-partitionable.
We next show how for some suitable sets $F$, the hypergraphs $Q_{pq}^{F}$ are
packing. To do so, we start with the following Lemma.

\begin{Lemma}\label{b_QFpq} Let $t$ be an element of the blocker
$b\left(  Q_{pq}^{F}\right)$.

\begin{description}
\item{\rm{I.}} If $F\subseteq\left(  F_{pq}\right)  ^{\ast}$ then
$t$ is of the form (\ref{eq_bQpq}.a-d)
or $t$ is of type (\ref{eq_bQpq}.e-h) whenever $t\nsubseteq f^{\ast}$, $f\in F$.

\item{\rm{II.}} If $F=\left(  E_{pq}\backslash F^{\prime}\right)  ^{\ast}$ where $F^{\prime}\subseteq F_{pq},F^{\prime}\neq\emptyset$ then $t$ is of the form (\ref{eq_bQpq}.a,b) or $t \in F_{pq}\backslash F^{\prime}$.
\end{description}
\end{Lemma}

\demo Since $E\left(  Q_{pq}^{F}\right)  \supseteq E_{pq}$, we have
that for every $t^{\prime}\in b\left(Q_{pq}^{F}\right)  $ $\exists
t\in b\left(  Q_{pq}\right)  $ such that $t^{\prime}=t\cup
x^{\prime}$. Let $t^{\prime}$ be an element in
$b\left(Q_{pq}^{F}\right)  $.

\medskip

Case $\rm{I}$. If $t \in b(Q_{pq})$ of type (\ref{eq_bQpq},a-d) then
$t \in b(Q_{pq}^{F})$ and $t^{\prime}=t$. Now, for $t=Abc$ a minimal
transversal of type (\ref{eq_bQpq}.e-h), if
$f=A^{\ast}B^{\ast}c^{\ast}\in F$, as $t^{\prime}\cap
f\neq\emptyset$ and $t\cap f=\emptyset$ then $x^{\prime}\cap
f=\emptyset$. If $x^{\prime}\cap A^{\ast}c^{\ast}\neq\emptyset$ then
$t^{\prime}$ contains $t \in b(Q_{pq})$ of the form
(\ref{eq_bQpq}.a). Thus $x^{\prime}\subseteq B^{\ast}$, but this
implies that $t^{\prime}$ contains a $t \in b(Q_{pq})$ of the form
(\ref{eq_bQpq}.a,c or d). Therefore $f\notin F$, and $t\in b\left(
Q_{pq}^{F}\right)$ implying $t^{\prime}=t$.

\medskip

Case $\rm{II}$. Notice that for any $t \in b(Q_{pq})$ of types
(\ref{eq_bQpq}.a,b) we have that $t \in b(Q_{pq}^{F})$. Let
$t=a_{i}a_{j}^{\ast}c$ be a minimal transversal of types
(\ref{eq_bQpq}.c,d), and choose $A_{1}\subset A$ such that $a_{i}\in
A_{1},a_{j}\notin A_{1}$. Since, $A_{1}^{c}A_{1}^{\ast}Bc^{\ast}$
and $A_{1}^{c}A_{1}^{\ast}B^{\ast}c^{\ast}$ are both edges of
$Q_{pq}^{F}$, we have that, $t^{\prime}$ could only be a minimal
transversal of $Q_{pq}^{F}$ if
$a_{i}^{\ast},a_{j},b_{k}b_{l}^{\ast}$ or $c^{\ast}$ are contained
in $x^{\prime}$, which implies that $t^{\prime}$ contains at least
one $t \in b(Q_{pq})$ of type (\ref{eq_bQpq}.a,b).

Now, let $t=Ab_{i}c$ be a minimal transversal of type
(\ref{eq_bQpq}.e-h). If $\emptyset\neq x^{\prime}\cap
A^{\ast}B^{\ast}c^{\ast}$, this implies that $t^{\prime}$ contains a
$t \in b(Q_{pq})$ of type (\ref{eq_bQpq}.a), or we are in the case
$t$ is a minimal transversal of type (\ref{eq_bQpq}.c,d). Hence
$x^{\prime}\subseteq B$, and $f=A^{\ast}B^{\ast}c^{\ast}\notin F$.

On the other hand, if $B_{1}$ is a proper subset of $B$ then
$A^{\ast}B_1\left(B_1^{c}\right )^{\ast}c^{\ast}$ is an edge of
$Q_{pq}^{F}$. But this implies that $ABc\subseteq t^{\prime}$,
therefore $t^{\prime}=ABc$, which is a minimal transversal of
$Q_{pq}^{F}$. \QED

\medskip

We construct a new class of clutters with the packing property
obtained by adding some new hyperedges to the hypergraphs $Q_{pq}$.

\newpage

\begin{Theorem}\label{QFpq_Packing} The hypergraphs $Q_{pq}^{F}$ have the packing property for:

\begin{description}
\item{\rm{I.}} $F\subseteq \left(F_{pq}\right)^{\ast}$, $F\neq\emptyset$, or

\item{\rm{II.}} $F=\left(  E_{pq}\backslash F^{\prime}\right)  ^{\ast}$ with $F^{\prime
}\subseteq F_{pq},F^{\prime}\neq\emptyset$.
\end{description}
\end{Theorem}

\demo From the construction of $Q_{pq}^{F}$ we have that
$\tau\left(Q_{pq}^{F}\right)  =\nu\left(  Q_{pq}^{F}\right)  =2$.
Now let $q_{pq}^{F}:=Q_{pq}^{F}\backslash W/Z$ be a proper minor of
$Q_{pq}^{F}$. If $\emptyset \in E(q_{pq}^{F})$ or $\emptyset =
E(q_{pq}^{F})$ then by convention it packs. Therefore, for both
cases we must show that for every proper minor $q_{pq}^{F}$ such
that $\emptyset \notin E(q_{pq}^{F})$ and $\emptyset \neq
E(q_{pq}^{F})$ the following equality holds:
\begin{equation}
1\leq \tau\left(  q_{pq}^{F}\right)  =\nu\left(  q_{pq}^{F}\right)  \label{eq_PP}%
\end{equation}

For the proof, we will use the fact that
$b\left(  H\backslash v\right)  =\left(  b\left(  H\right)  /v\right)  ^{\min}$,
$b\left(H/v\right)  =b\left(  H\right)  \backslash v$ and that deletion and
contraction are associative and commute. We will do first the contractions,
and then deletions. If $\tau\left(  q_{pq}^{F}\right)=1$, as
$\emptyset \notin E(q_{pq}^{F})$ and $\emptyset \neq E(q_{pq}^{F})$ then
$q_{pq}^{F}$ packs. Thus we can assume $\tau\left(  q_{pq}^{F}\right)\geq 2$.

We have that
if $t \in b\left(Q_{pq}^{F}\right)$, then there exists $t^{\prime}\in b\left(Q_{pq}\right)$ such that
$t^{\prime}\subseteq t$,
therefore,
$\tau\left(  Q_{pq}\backslash W/Z\right)  \leq\tau\left(q_{pq}^{F}\right)  $.

If $\tau\left(  q_{pq}^{F}\right)  =\tau\left(Q_{pq}\backslash
W/Z\right)  $ then $\nu\left(  q_{pq}^{F}\right)=\tau\left(
q_{pq}^{F}\right)  $ since $E\left(  q_{pq}^{F}\right)  \supseteq
E\left(  Q_{pq}\backslash W/Z\right)  $ and $Q_{pq}$ packs. Hence,
we must prove that $\tau\left(  q_{pq}^{F}\right)  =\nu\left(
q_{pq}^{F}\right)$ when
$\tau\left(q_{pq}^{F}\right)>\tau\left(Q_{pq}\backslash W/Z\right)$.

%\medskip
%\begin{enumerate}{}
%
%\item [\bf Case I] ($F\subseteq \left(F_{pq}\right)^{\ast}$ with $F\neq\emptyset$):

\medskip

\noindent
{\bf [Case I]} ($F\subseteq \left(F_{pq}\right)^{\ast}$ with $F\neq\emptyset$):
Notice that $\tau\left(q_{pq}^{F}\right)>\tau\left(Q_{pq}\backslash W/Z\right)$
could happen only if
$\tau\left(  Q_{pq}\backslash W/Z\right)=\left\vert Abc\backslash W\right\vert $ and
$Z\cap Abc=\emptyset$, where
$Abc$ is a minimal transversal of type (\ref{eq_bQpq}.e-h),
and $Abc\notin b\left(  q_{pq}^{F}\right)$. But this implies that
$A^{\ast}B^{\ast}c^{\ast}\in F$.

\medskip

%\begin{enumerate}
%
%\item [Case I.1]

%\begin{description}
%
%\item
\noindent
{\bf [I.1]} If $Z\cap vv^{\ast}\neq\emptyset$ for all $v\in PQr$, then
$A_{1}\left(A_{1}^{c}\right)  ^{\ast}B_{1}\left(  B_{1}^{c}\right)^{\ast}d\subseteq Z$
where $A_{1}\subseteq A,B_{1}\subseteq B$ and $d\in rr^{\ast}$.
As $Z\cap Abc=\emptyset$, then $d=c^{\ast}$, $A_{1}=\emptyset$ and $b\notin B_{1}$.
Therefore $Z=A^{\ast}B_{2}B_{3}^{\ast}c^{\ast}$, where
$B_{1}\subseteq B_{2}\subseteq B$,
$\left(  B_{1}^{c}\right)  ^{\ast}\subseteq B_{3}^{\ast}\subseteq B^{\ast}$ and
$b\notin B_{2}$. If $B^{\ast}=B_{3}^{\ast}$ then
$\emptyset\in E\left(  q_{pq}^{F}\right)$, but we assumed that
$\emptyset\notin E\left(  q_{pq}^{F}\right)$.
Thus, there exists
$b_1^{\ast}\in B^{\ast}$ such that $b_1^{\ast}\notin B_{3}^{\ast}$, and consequently
$B_{2}\neq\emptyset$.

\medskip
%\item
%\noindent

By the form of the elements of $b\left(  Q_{pq}\right)$ and
of $b\left(  Q_{pq}^{F}\right)$ we have that:\newline
$b\left(  Q_{pq}/Z\right)=\left\{  b_{i}b_{j}^{\ast}c,Abc\right\} \backslash B_{2}B_{3}^{\ast}$,
$b\left(  Q_{pq}^{F}/Z\right)  =\left\{  b_{i}b_{j}^{\ast}c\right\}
\backslash B_{2}B_{3}^{\ast}$ and
$bb_{1}^{\ast}c\in b\left(  Q_{pq}/Z\right)
\cap b\left(  Q_{pq}^{F}/Z\right)$.
Therefore
$\left\vert bb_{2}^{\ast}c\backslash W\right\vert >\left\vert Abc\backslash W\right\vert$ for all
$b_{2}^{\ast}\notin B_{3}^{\ast}$. This is only possible if $A\subseteq W$,
and $b_{2}^{\ast}\notin W$ for all $b_{2}^{\ast}\notin B_{3}^{\ast}$;
\emph{i.e.} $W\cap\left(  B_{3}^{\ast}\right)  ^{c}=\emptyset$.
Consequently,
$3\geq\tau\left(  q_{pq}^{F}\right)  =\tau\left(  Q_{pq}\backslash W/Z\right)
+1\geq1$. Thus, either $W=AB_{4}$ or $W=AB_{4}c$, where
$B_{4}\subseteq B\backslash B_{2}$.
Since $B_{2}\neq\emptyset$ then $B_{4}\neq B$. Even more,
$\emptyset\neq f=\left(  B_{3}^{c}\right)  ^{\ast} = A^{\ast}B^{\ast}c\backslash
A^{\ast}B_{2}B_{3}^{\ast}c^{\ast}\in E\left(  q_{pq}^{F}\right) \backslash E\left(  Q_{pq}\backslash W/Z\right)$. Now, by the form of the
elements of $E\left(  Q_{pq}\right)$ we have that
$E\left(  Q_{pq}\backslash W/Z\right) =\left(  \left\{  V\left(  V^{c}\right)  ^{\ast}c:V\subset B\right\}  \cup B\right)  /B_{2}B_{3}^{\ast}\backslash W^{\prime}$
where
$W^{\prime}=W\backslash A$ and
$E\left(  q_{pq}^{F}\right) \supseteq E\left(  Q_{pq}\backslash W/Z\right)  \cup\left(  B_{3}^{c}\right)  ^{\ast}$.

%\end{description}

%\begin{enumerate}
%
%\item [subcase I.1.1]
%\noindent
\medskip

\noindent
{\bf [I.1.1]} If $c\in W^{\prime}$ then
$E\left(  Q_{pq}\backslash W/Z\right)  =\left\{ B\right\}  /B_{2}\backslash B_{4}$. Now, $\left(  B_{3}^{c}\right)  ^{\ast}$ is
independent with the edges of $Q_{pq}\backslash W/Z$, which implies that $q_{pq}^{F}$
packs.

\medskip
%\item [subcase I.1.2]
\noindent
{\bf [I.1.2]} Now if $c\notin W^{\prime}$ then
$E\left(  Q_{pq}\backslash W/Z\right) \supseteq\left(  \left\{  V\left(  V^{c}\right)  ^{\ast}c:V\subseteq B_{4}^{c}\right\}  \right)  /B_{2}B_{3}^{\ast}$;
since $B_{2}\cap B_{4}=\emptyset$, we have $B_{2}\left(  B_{2}^{c}\right)  ^{\ast}c\setminus B_{2}B_{3}^{\ast}\in E\left(
Q_{pq}\backslash W/Z\right)  $, but $B_{1}\subseteq B_{2}$, therefore
$B_{2}^{c}\subseteq B_{1}^{c}\subseteq B_{3}$.
Thus $c=B_{2}\left(  B_{2}%
^{c}\right)  ^{\ast}c\backslash B_{2}B_{3}^{\ast}\in E\left(  Q_{pq}\backslash
W/Z\right)  $. We have the following two remaining cases.

\medskip
%\begin{enumerate}
%
%\item [subcase I.1.2.1]
%\noindent

If $B_{4}\neq\emptyset$, $E\left(  Q_{pq}\backslash W/Z\right)
=\left(  \left\{  V\left(  V^{c}\right)  c:V\subseteq B_{4}^{c}\right\}
\right)  /B_{2}B_{3}^{\ast}$, then $\tau\left(  Q_{pq}\backslash W/Z\right)
=1$ and $c$ is independent with $\left(  B_{3}^{c}\right)  ^{\ast}$ hence,
$q_{pq}^{F}$ packs.

\medskip
%\item [subcase I.1.2.2]
%\noindent

If $B_{4}=\emptyset$, $B_{2}^{c}=\left\{  B\right\}
/B_{2}B_{3}^{\ast}\in E\left(  Q_{pq}\backslash W/Z\right)  $ then $B_{2}%
^{c},c,\left(  B_{3}^{c}\right)  ^{\ast}$ are independent in $q_{pq}^{F}$. As
$\tau\left(  q_{pq}^{F}\right)  \leq 3$ then $q_{pq}^{F}$ packs.

%\end{enumerate}
%\end{enumerate}
%\end{enumerate}

\medskip

%\begin{enumerate}

%\item [Case I.2]

%\begin{description}

%\item
\noindent
{\bf [I.2]} If $Z\cap vv^{\ast}=\emptyset$ for some $v\in PQr$ then, $2\geq\left\vert
vv^{\ast}\backslash W\right\vert \geq\tau\left(  q_{pq}^{F}\right)
>\tau\left(  Q_{pq}\backslash W/Z\right)  =\left\vert Abc\backslash
W\right\vert $ for all $v$ such that $Z\cap vv^{\ast}=\emptyset$. If
$Abc\subseteq W$, as $Abc\in b\left(  Q_{pq}\right)  $ then $E\left(
q_{pq}^{F}\right)  =\left\{  A^{\ast}B^{\ast}c^{\ast}\right\}  /Z\backslash
W$, and $q_{pq}^{F}$ packs. Assume $Abc\nsubseteq W$, then $\tau\left(
Q_{pq}\backslash W/Z\right)  =1$ and $\tau\left(  q_{pq}^{F}\right)  =2$. Notice that if $u\in W$, then $u^{\ast}\in Z$; since otherwise, knowing that
$W\cap Z = \emptyset$ then $uu^{\ast}\cap Z=\emptyset$. Thus $1\geq\left\vert
uu^{\ast}\backslash W\right\vert >\left\vert Abc/W\right\vert =\tau\left(
Q_{pq}\backslash W/Z\right)  $ which is a contradiction.

%\end{description}

\medskip

%\begin{enumerate}
%
%
%\item [subcase I.2.1]
\noindent
{\bf [I.2.1]} If $bc \subseteq W$, then
$b_{i}^{\ast}\in b\left(  Q_{pq}{F}\backslash W \right )$. Therefore
$b_{i}^{\ast}\in B^{\ast}\subset Z$. In this case $Abc\backslash W =\{a\}\subseteq A$, and $c^{\ast}\in Z $; implying
$\left(  A\backslash\left\{a\right\}  \right)  ^{\ast}B^{\ast}c^{\ast}\subseteq Z$ and
$\left( A\backslash\left\{  a\right\}  \right)  bc\subseteq W$. Now, by construction of
$E\left(  Q_{pq}\right)  $, we have that $E\left(  q_{pq}^{F}\right)
=\left\{  a,a^{\ast}\right\}  \backslash\left(  W\cap Ba^{\ast}\right)
/\left(  Z\cap Ba\right)  $, but $\tau\left(  q_{pq}^{F}\right)  =2$, and
consequently $q_{pq}^{F}$ packs.

\medskip

%\item [subcase I.2.2]
\noindent
{\bf [I.2.2]} Now, if $bc\nsubseteq W$, then $A\subseteq W$
and $A^{\ast}\subseteq Z$. Even more, $b\in W$ or $c\in W$.

%\begin{enumerate}
%\item [subcase I.2.2.1]
%\noindent

\medskip

If $c\in W$, then
$c^{\ast}\in Z$. Therefore $Ac\subseteq W$ and $A^{\ast}c^{\ast}\subseteq Z$. By the form of the elements of $E\left(  Q_{pq}\right)$ and of $F$ we have that
$E\left(  q_{pq}^{F}\right) =\left\{  B,B^{\ast}\right\} \backslash\left(
W\cap BB^{\ast}\right)/\left(  Z\cap BB^{\ast}\right)$, hence $q_{pq}^{F}$
packs.

%\item [subcase I.2.2.2]
%\noindent

\medskip

If $b\in W$, then $Ab\subseteq W$ and
$A^{\ast}b^{\ast}\subseteq
Z$. So the edge set of $q_{pq}^{F}$ is given by $E\left(  q_{pq}^{F}\right)  =
\left(  \left\{  V\left(  V^{c}\right)
^{\ast}c:V\subseteq B\backslash b\right\}  \cup\left(  B\backslash b\right)
^{\ast}c^{\ast}\right)  \backslash\left(  W\cap BB^{\ast}c^{\ast}\right)  /Z$.
Since
$E\left(  q_{pq}^{F}\right)  \neq E\left(  Q_{pq}\backslash W/Z\right)$ we have that $W\cap\left(  B\backslash b\right)  ^{\ast}c^{\ast} = \emptyset$. Therefore $W=AB_{5}$, $Z\supseteq A^{\ast}B_{6}^{\ast}$ with
$b\in B_{5}\subseteq B_{6}\subseteq B$. As $\tau\left(  Q_{pq}\backslash
W/Z\right)  =1$, we must have $B_{5}^{c}c\backslash Z=\left(  B_{5}%
^{c}\right)  B_{5}^{\ast}c\backslash Z\in E\left(  Q_{pq}\backslash
W/Z\right)  $. Then $B_{5}^{c}c\backslash Z$ and $\left(  B\backslash b\right)  ^{\ast
}c^{\ast} \backslash Z$ are independent in $E\left(  q_{pq}^{F}\right)  $ and hence,
$q_{pq}^{F}$ packs.

%\end{enumerate}
%\end{enumerate}
%\end{enumerate}

%\begin{enumerate}

%\item [\bf Case II] ( $F=\left(  E_{pq}\backslash F^{\prime}\right)  ^{\ast}$ with $F^{\prime
%}\subseteq F_{pq},F^{\prime}\neq\emptyset$):

\bigskip

\noindent {\bf [Case II]} ( $F=\left(  E_{pq}\backslash F^{\prime}\right)  ^{\ast}$ with $F^{\prime
}\subseteq F_{pq},F^{\prime}\neq\emptyset$):

\medskip

\noindent{\bf [II.1]} If $\tau\left(  Q_{pq}^{F}/Z\right)  =4$, then
$vv^{\ast}\cap Z\neq\emptyset$ for all $v\in PQr$. Thus
$Z=P_{1}P_{2}^{\ast}Q_{1}Q_{2}^{\ast }D$, where $P_{1}\cup P_{2}=P$,
$Q_{1}\cup Q_{2}=Q$ and $\emptyset\neq D\subseteq rr^{\ast}$. Since
$4=\tau\left(  Q_{pq}^{F}/Z\right)  =\left\vert
p_{1}p_{2}^{\ast}q_{1}q_{2}^{\ast}\right\vert $, $p_{1}p_{2}^{\ast}q_{1}%
q_{2}^{\ast}\cap Z=\emptyset$ then $P_{1},P_{2}\neq P$ and $Q_{1},Q_{2}\neq
Q$. Moreover, $P_{1},P_{2}\neq\emptyset$, $Q_{1},Q_{2}\neq\emptyset$. By the
form of $E\left(  Q_{pq}^{F}\right)  $ we have that $E\left(  \left(
Q_{pq}^{F}/Z\right)  ^{\min}\right)  =\left\{  P_{1}^{c},\left(  P_{2}^{\ast
}\right)  ^{c},Q_{1}^{c},\left(  Q_{2}^{\ast}\right)  ^{c}\right\}  $.
Therefore $q_{pq}^{F}$ packs.

\medskip

\noindent{\bf [II.2]} If $\tau\left(  Q_{pq}^{F}/Z\right)  =p+q+1$,
then $\tau\left( Q_{pq}^{F}/Z\right)  =\left\vert ABc\right\vert $
and $ABc\cap Z=\emptyset$, implying $A^{\ast}B^{\ast}c^{\ast}\notin
F$. Moreover, since $vv^{\ast}\cap Z\neq\emptyset$ for all $v\in
PQr$, then $Z=A^{\ast}B^{\ast}c^{\ast}$. By the form of $E\left(
Q_{pq}^{F}\right)  $ we have that $E\left(  \left(
Q_{pq}^{F}/Z\right)  ^{\min}\right)  =\left\{  \left\{  v\right\}
:v\in ABc\right\}  $. Therefore $q_{pq}^{F}$ packs.

\medskip

\noindent{\bf [II.3]} If $\tau\left(  Q_{pq}^{F}/Z\right)  =2$, then
$\tau\left( Q_{pq}^{F}/Z\right)  =\left\vert vv^{\ast}\right\vert $
for some $v\in PQr$; \emph{i.e.} $vv^{\ast}\cap Z=\emptyset$. Let us
consider the case where
\begin{equation}
2=\tau\left(  q_{pq}^{F}\right)  =\left\vert vv^{\ast}\backslash
W\right\vert
>\tau\left(  Q_{pq}\backslash W/Z\right)  \label{eq_case3}%
\end{equation}

By Proposition~\ref{2-part} we have that $W\cap uu^{\ast}=\emptyset$
for all $u$ such that $uu^{\ast}\cap Z=\emptyset$. So, if $u\in W$,
then $u^{\ast}\in Z$ since otherwise $\tau\left(  q_{pq}^{F}\right)
<2$, a contradiction.

For the proof, let $W=A_{1}A_{2}^{\ast}B_{1}B_{2}^{\ast}D_{1}$ and
$Z=A_{3}A_{4}^{\ast}B_{3}B_{4}^{\ast}D_{2}$, with $A_{1}\subseteq
A_{4}$, $A_{2}\subseteq A_{3}$, $B_{1}\subseteq B_{4}$,
$B_{2}\subseteq B_{3}$ and $D_{1}^{\ast}\subseteq D_{2}\subseteq
rr^{\ast}$. Furthermore, since $W\cap Z=\emptyset$ we have that
$A_{1}\cap A_{3}=\emptyset$, $A_{2}\cap A_{4}=\emptyset$, $B_{1}\cap
B_{3}=\emptyset$, $B_{2}\cap B_{4}=\emptyset$ and $D_{1}\cap
D_{2}=\emptyset$. Moreover, $A_{1}\cap A_{2}=\emptyset$, $B_{1}\cap
B_{2}=\emptyset$ and $\left\vert D_{1}\right\vert \leq1$.

On the other hand, some of the sets $A_{1}$, $A_{2}$, $B_{1}$ or
$B_{2}$ is empty. This is necessary, since otherwise this would
imply that $\emptyset\in b\left(  q_{pq}^{F}\right)  $, which is a
contradiction. Therefore, assume without loss of generality that
$B_{2}=\emptyset$. For the rest of the proof assume that
$f_{1},f_{2}\in E\left(  q_{pq}^{F}\right)  $.

\medskip

\noindent{\bf [II.3.1]} Assume that  $A_{1}\neq\emptyset$ and
$A_{2}\neq\emptyset$, then:

\medskip

\noindent{\bf [II.3.1.1]} If $B_{1}\neq\emptyset$, we have
$B_{4}^{\ast}=B^{\ast}$, hence:

\medskip

If $D_{1}=\emptyset$; by the form of $E\left( Q_{pq}^{F}\right) $,
there exists $f_{1}\subseteq\left( A_{4}^{c}\backslash A_{2}\right)
^{\ast}c$ and $f_{2}\subseteq\left( A_{3}^{c}\backslash A_{1}\right)
c^{\ast}$ then $E\left( q_{pq}^{F}\right)  $ packs.

\medskip

If $D_{1}\neq\emptyset$; assume $D_{1}=c$. If $A_{1}\cup A_{2}=A$
then $\emptyset\in E\left(  q_{pq}^{F}\right) $. Thus, suppose that
$A_{1}\cup A_{2}\neq A$. By the form of the elements of $E\left(  Q_{pq}%
^{F}\right)  $, there exists $f_{1}\subseteq\left(
A_{4}^{c}\backslash A_{2}\right)  ^{\ast}$ and $f_{2}\subseteq\left(
A_{3}^{c}\backslash A_{1}\right)  $, therefore $E\left(
q_{pq}^{F}\right)  $ packs.

\medskip

\noindent{\bf [II.3.1.2]} If $B_{1}=\emptyset$. Consider the case
where
$D_{1}=\emptyset$, then there are $f_{1},f_{2}$ such that $f_{1}%
\subseteq\left(  A_{3}^{c}\backslash A_{1}\right)  Bc$ and $f_{2}%
\subseteq\left(  A_{4}^{c}\backslash A_{2}\right)
^{\ast}B^{\ast}c^{\ast}$. Then $q_{pq}^{F}$ packs.

\medskip

On the other hand, if $D_{1}\neq\emptyset$, assume $D_{1}=c$. This
implies that $c^{\ast}\in Z$. For this case there are $f_{1},f_{2}$
such that $f_{1}\subseteq\left(  A_{3}^{c}\backslash A_{1}\right) B$
and $f_{2}\subseteq\left(  A_{4}^{c}\backslash A_{2}\right)
^{\ast}B^{\ast}$.

\medskip

\noindent{\bf [II.3.2]} Assume that $A_{1}=\emptyset$ or
$A_{2}=\emptyset$. Without loss of generality let $A_{2}=\emptyset$, and let $A_{1},B_{1}%
\neq\emptyset$.

\medskip

\noindent{\bf [II.3.2.1]} So if $A_{1}\neq A$ or $B_{1}\neq B$ then:

If $D=\emptyset$, then there exists $f_{1}\subseteq A_{1}^{c}\left(
B_{4}^{c}\right) ^{\ast}c$ and $f_{2}\subseteq\left(
A_{4}^{c}\right) ^{\ast}B_{1}^{c}c^{\ast}$, then $E_{pq}^{F}$ packs.

\medskip

If $D\neq\emptyset$, let $D=c$,
implying $c^{\ast}\in Z$. Then there are $f_{1},f_{2}\in E\left(  q_{pq}%
^{F}\right)  $ such that $f_{1}\subseteq A_{1}^{c}\left(
B_{4}^{c}\right) ^{\ast}$ and $f_{2}\subseteq\left(
A_{4}^{c}\right)  ^{\ast}B_{1}^{c}$.

\medskip

\noindent{\bf [II.3.2.2]} So if $A_{1}=A$ and $B_{1}=B$, then
$D=\emptyset$ and $f_{1}=c$ and $f_{2}=c^{\ast}$.

\medskip

\noindent{\bf [II.3.3]} Assume that $B_{1}=\emptyset$, $A_{2}=\emptyset$ and $A_{1}%
\neq\emptyset$, then:

\medskip

\noindent{\bf [II.3.3.1]} If $D_{1}=\emptyset$, then $f_{1}\subseteq A_{1}%
^{c}B^{\ast}c$ and $f_{2}\subseteq\left(  A_{1}^{c}\right)
^{\ast}Bc^{\ast}$.

\medskip

\noindent{\bf [II.3.3.2]} If $D_{1}\neq\emptyset$, then assume that
$D_{1}=c$ and $c^{\ast}\in Z$. If $A_{1}\neq A$, there are
$f_{1},f_{2}$ such that $f_{1}\subseteq\left(  B^{\prime}\right)
^{c}\left( B^{\prime}\right) ^{\ast}$ and $f_{2}\subseteq\left(
B^{\prime}\right)  \left(  \left( B^{\prime}\right)  ^{c}\right)
^{\ast}$ for some $\emptyset\neq B^{\prime }\subset B$, therefore
$q_{pq}^{F}$ packs.

\medskip

\noindent{\bf [II.3.4]} Assume that $B_{1}=\emptyset$, $A_{2}=\emptyset$ and $A_{1}%
\neq\emptyset$. Now let $W\supseteq rr^{\ast}$ then notice that for
this case there are $f_{1},f_{2}\in E\left(  Q_{pq}\backslash
W\right)  $ two independent edges, therefore, $f_{1}\backslash Z$
and $f_{2}\backslash Z$ should contain independent edges in
$q_{pq}^{F}$, therefore $q_{pq}^{F}$ packs.

%Whenever $Z$ is such that
%$\tau\left(  Q_{pq}^{F}/Z\right)  =2$, we have that there must be $f_{1}%
%,f_{2}\in E\left(  Q_{pq}^{F}/Z\right)  $ such that $f_{i}\subseteq
%e_{i},i=1,2$, just like in the first case, implying once more that $q_{pq}%
%^{F}$ is packing.
%
%\medskip
%
%Now, if $\tau\left(  Q_{pq}^{F}/Z\right)  =4$, then
%$Z=\mathbb{PQ}rr^{\ast}$ where $\mathbb{P}=P_{1}\left(  P_{1}^{c}\right)
%^{\ast}$ and $\mathbb{Q}=Q_{1}\left(  Q_{1}^{c}\right)  ^{\ast}$ for some pair
%$P_{1},Q_{1}$ proper subsets of $P$ and $Q$ respectively.
%
%%\medskip
%%
%%But $Q_{pq}^{F}$
%%coincides with ..., and therefore $q_{pq}^{F}$ is packing.
%
%
%\medskip
%
%Finally, let us assume that $\left\vert ABc\right\vert =\tau\left(  Q_{pq}%
%^{F}/Z\right)  $, this implies that $Z=\left(  ABc\right)  ^{\ast}$.
%
%\medskip
%
%Consequently, $Q_{pq}^{F}/Z=\left\{  v:v\in ABc\right\}  $. In other words we
%have an identity matrix of order $p+q+1$, which is a packing hypergraph.

\bigskip

Concluding with this the proof that $Q_{pq}^{F}$ has the packing
property. \QED

%\end{enumerate}
%\end{enumerate}

\medskip

As mentioned before, by construction all the clutters in the class $Q_{pq}^{F}$ are $2$-partitionable.
In particular for the clutters given in Theorem~\ref{QFpq_Packing} we have:

\begin{Proposition}\label{Deltar-Q6}
Consider $Q_{pq}^{F}$ for $F$ as in cases $I$ and $II$ of Theorem~\ref{QFpq_Packing}. Then, $\Delta_r(B(Q_{pq}^{F}))=1$.
\end{Proposition}

\demo We have that $Q_{pq}^{F}$ is 2-partitionable and has the
K\"onig property. By the observation before
Proposition~\ref{jan17-04} and Remark~\ref{RR}, we have that
$rank(A)=rank\left(B(Q_{pq}^{F})\right) \leq d+1$. In this case
$d=p+q+1$, so $rank\left(B(Q_{pq}^{F})\right) \leq p+q+2$.

On the other hand,

\medskip

\begin{center}
\begin{tabular*}{1.3in}
[r]{cp{0.25in}p{0.33in}p{0.2in}p{0.2in}}
&       $\quad P^{\ast}$  &   $\,\, Q^{\ast}$ & \hglue-0.06in $r^{\ast}$ & \hglue-0.03in $r$\\
\multicolumn{1}{l}{$L=$} & \multicolumn{4}{l}{$\left[
\begin{tabular*}{1.3in}
[c]{p{0.13in}p{0.2in}p{0.2in}p{0.2in}}%
$\mathbf{0}$ & $I$ & $\mathbf{1}$ & $\mathbf{0}$\\
$\mathbf{0}$ & $\mathbf{0}$ & $0$ & $1$\\
$I$ & $\mathbf{0}$ & $\mathbf{0}$ & $\mathbf{1}$\\
$\mathbf{1}$ & $\mathbf{0}$ & $1$ & $0$%
\end{tabular*}
\right]  $}%
\end{tabular*}
\end{center}

\medskip

\noindent is a submatrix of $B(Q_{pq}^{F})$, and $L$ is reducible by elementary matrix transformations to $I_{p+q+2}$ (the identity matrix
of order $p+q+2$). Hence ${\rm det}(L)=1$ and it follows that
$\Delta_r(B(Q_{pq}^{F}))=1$.
\QED

\medskip

Note that the previous result gives support to Conjecture~\ref{git-val-vi}.

\section{A new infinite family of mengerian clutters}

In this section we will prove that if $F=F_{pq}^{\ast}$ then $Q_{pq}^{F}$ is a mengerian hypergraph.

Among the general classes of clutters known to verify the Conjecture
of Conforti and Cornu\'ejols are: \emph{binary}, \emph{balanced} and
\emph{dyadic} clutters. We now prove that the family $Q_{pq}^{F}$ of
Theorem~\ref{QFpq_Packing} does not belong to any of these classes.

\medskip

Let us denote by $\bigtriangleup$ the symmetric difference operator. A
hypergraph $H=\left(  V,E\right)$ is:

\begin{description}

\item{a)} \emph{Binary} if for every
$e_{1},e_{2},e_{3}\in E$ there is an $e\in E$ such that $e\subseteq
e_{1}\bigtriangleup e_{2}\bigtriangleup e_{3}$.

\item{b)} \emph{Dyadic} if for every pair $\left(  e,t\right)  $ with $e\in E,t\in
b\left(  H\right)  $ the inequality $\left\vert t\cap e\right\vert \leq2$
holds.

\item{c)} \emph{Balanced} if no square submatrix
of odd order contains exactly two $1^{\prime}s$ per row and per column.

\end{description}

\begin{Lemma}\label{notdbb}
The hypergraphs $Q_{pq}^{F}$ of Theorem~\ref{QFpq_Packing} are not dyadic, nor
binary, nor balanced.
\end{Lemma}

\demo For the proof let $e,f\in E_{pq}$ such that $e^{c}\in
F,f^{c}\notin F$ (the existence of both elements is a consequence of
the construction given in Theorem~\ref{QFpq_Packing}). Now, to prove
that $Q_{pq}^{F}$ is not dyadic, notice that there is an edge $e\in
E_{pq}$ and a minimal transversal $t$ of type (\ref{eq_bQpq}.b) such
that $\left\vert t\cap f\right\vert \geq3$. Moreover, notice that
$e\bigtriangleup e^{c}\bigtriangleup f=f^{c}\notin E\left(
Q_{pq}^{F}\right)  $ so $Q_{pq}^{F}$ is not binary. Now, to prove
that $Q_{pq}^{F}$ is not balanced, notice that the edges $PQr,$
$PQ^{\ast}r^{\ast}$ and $P^{\ast}Qr^{\ast}$ intersected with $p\in
P,$ $q\in Q$ and $r^{\ast}$ constitute a counterexample. \QED

\medskip

Observe that if we consider for $Q_{pq}^{F}$ that $F=E_{pq}^{c}$ (which was
not taken into account for the second case of Theorem~\ref{QFpq_Packing}), then $Q_{pq}^{F}$
is a binary hypergraph.

\medskip

\begin{Theorem}\label{menger}
Let $H$ be an ideal hypergraph. Then the following statements are equivalent.

\begin{enumerate}
\item $H$ is mengerian

\item If $\tau^{w}\left(  H\right)  >0$, then there exists $e\in E\left(
H\right)  $ such that $\tau^{w}\left(  H\right)  =\tau^{w-e}\left(  H\right)
+1$.
\end{enumerate}
\end{Theorem}

\demo $\left(  2\Rightarrow1\right)  $ As $H$ is ideal,
$\tau^{w}\left(  H\right) =wx$, where $x\in b\left(  H\right)  $. We
only need to prove that $\tau ^{w}\left(  H\right)  =\nu^{w}\left(
H\right)  =y\mathbf{1}$, where $Ay\leq w$, $y\geq0$ and $y$ integer.
We will prove this by induction over $\tau ^{w}\left(  H\right)  $.
If $\tau^{w}\left(  H\right)  =0$, then $y=0$ satisfies the
proposition. Now suppose this true for $\tau^{w}\left( H\right) <n$.
If $\tau^{w}\left(  H\right)  =n$ by $2)$ we know that
$\tau^{w}\left(  H\right)  =\tau^{w-e}\left(  H\right)  +1$, where
$e\in E\left(  H\right)  $. Suppose that $e=v_{k}$, then by
induction hypothesis there exists $y^{\prime}=\left(
y_{1}^{\prime},\ldots,y_{k}^{\prime}\right) $ such that
$\tau^{w-e}\left(  H\right)  =y^{\prime}\mathbf{1}$ and
$Ay^{\prime}\leq w-e$, $y^{\prime}\geq0$, $y_{1}^{\prime}v_{1} +
\cdots + y_{k}^{\prime}v_{k} \leq w-v_{k}$. Let $y=\left(
y_{1}^{\prime},\ldots,y_{k}^{\prime}+1\right)  $ then $Ay\leq w$,
$y\geq0$ and integer and
$y\mathbf{1}=y^{\prime}\mathbf{1}+1=\tau^{w-e}\left(  H\right)
+1=\tau^{w}\left(  H\right)  $. But this implies $H$ is mengerian.

$\left(  1\Rightarrow2\right)  $ As $H$ is mengerian,
$\tau^{w}\left( H\right)  =wx=\nu^{w}\left(  H\right)  =y\mathbf{1}$
where $x\geq0$, $xA\geq\mathbf{1}$, $y\geq0$, $Ay\leq w$ and $x$,
$y$ integer vectors. Suppose that $\tau^{w}\left(  H\right)  >\ 0$.
If $y=\left(  y_{1},\ldots ,y_{k}\right)  $, we have that
$y_{1}v_{1}+\cdots+y_{k}v_{k}\leq w$ and as $0<\tau^{w}\left(
H\right)  =\nu^{w}\left(  H\right)  =y_{1}+\cdots+y_{k}$, then
$y_{i}>0$ for some $i$. Suppose $y_{k}$ is precisely that element.
Then define $w^{\prime}=w-v_{k}>0$, and $y^{\prime}=\left(
y_{1},\ldots
,y_{k}-1\right)  $. Then $y^{\prime}\geq0$, $Ay^{\prime}\leq w-v_{k}%
=w^{\prime}$, and $\nu^{w^{\prime}}\left(  H\right)  \geq y^{\prime}%
\mathbf{1}=y\mathbf{1}-1$. Thus $\nu^{w^{\prime}}\left(  H\right)
\geq\nu ^{w}\left(  H\right)  -1$ and $\tau^{w^{\prime}}\left(
H\right)  \geq\tau ^{w}\left(  H\ \right)  -1$. Therefore
$\tau^{w}\left(  H\right)  =wx=\left( w^{\prime}+v_{k}\right)
x=w^{\prime}x+v_{k}x\geq\tau^{w^{\prime}}\left( H\right)
+v_{k}x\geq\tau^{w}\left(  H\right)  +v_{k}x-1$. This implies,
$1\geq v_{k}x$, but $xA\geq \mathbf{1}$ consequently $v_{k}x=1$ and
$\tau^{w}\left( H\right)  =\tau^{w^{\prime}}\left(  H\right) +1$.
\QED

\medskip

Now, let $p_{\min}:=\min_{p\in P}\left\{  w\left(  p\right) \right\}
$ and define $p_{\min}^{\ast}$, $q_{\min}$ and $q_{\min}^{\ast}$ in
an analogous way. We construct the sets $P_{>} =\left\{p\in
P:w\left( p\right)>p_{\min}\right\}  $ and $P_{>}^{\ast} = \left\{
p^{\ast} \in P^{\ast}:w\left( p^{\ast} \right)>p_{\min}^{\ast}
\right\} $. In an analogous manner we construct the sets  $Q_{>}$
and $Q_{>}^{\ast}$. Now take $\mathbb{P\subseteq }P_{>}\cup
P_{>}^{\ast}$ ($\mathbb{Q\subseteq}Q_{>}\cup Q_{>}^{\ast}$) with
maximum cardinality such that either $p_{i}$ or $p_{i}^{\ast}$
($q_{i}$ or $q_{i}^{\ast}$) but not both is an element of
$\mathbb{P}$ ($\mathbb{Q}$). These subsets fulfill the following
conditions.

\begin{Lemma}\label{Rk_2} If $p_{\min}=p_{\min}^{\ast}=0$ ($q_{\min}=q_{\min}^{\ast}=0$),
and $\tau^{w}>0$ then $p=\left\vert \mathbb{P}\right\vert $ ($q=\left\vert
\mathbb{Q}\right\vert $) and $q_{\min}+q_{\min}^{\ast}$, $w\left(  r\right)
\geq\tau^{w}$ ($p_{\min}+p_{\min}^{\ast}$, $w\left(  r^{\ast}\right)  \geq
\tau^{w}$).
\end{Lemma}

\demo If $p_{\min}=p_{\min}^{\ast}=0$ then $p=\left\vert
\mathbb{P}\right\vert $ (if not, there would be $p_{i},p_{i}^{\ast}$
such that $w\left(  p_{i}p_{i}^{\ast }\right)  =0$, which is a
contradiction). Now, that $q_{\min}+q_{\min}^{\ast}$ and $w\left(
r\right)  \ $are greater or equal to $\tau^{w}$ is a consequence of
(\ref{eq_bQpq}.b) and (\ref{eq_bQpq}.c). \QED

\medskip

In the sequel we denote by $\overline{t}$ a minimum weight $t \in b(Q_{pq}^F)$.

\medskip

\begin{Corollary}\label{Cor_2} For the case stated in the former Lemma, we can pick $e$ from:%
\begin{equation}%
\begin{array}
[c]{cc}%
e=\mathbb{P}Qr & \text{if }\left\vert \mathbb{P}\right\vert =p\text{ and
}q_{\min}>0\\
e=\mathbb{P}Q^{\ast}r & \text{if }\left\vert \mathbb{P}\right\vert =p\text{
and }q_{\min}^{\ast}>0\\
e=P\mathbb{Q}r^{\ast} & \text{if }\left\vert \mathbb{Q}\right\vert =q\text{
and }p_{\min}>0\\
e=P^{\ast}\mathbb{Q}r^{\ast} & \text{if }\left\vert \mathbb{Q}\right\vert
=q\text{ and }p_{\min}^{\ast}>0
\end{array}
\label{eq_e_Cor2}%
\end{equation}
such that $\tau^{w}=\tau^{w^{\prime}}+1$ and $w^{\prime}\geq0$.
\end{Corollary}

\demo That $w^{\prime}\geq0$ follows directly from the selection of
$e$. Thus, to prove that $\tau^{w}=\tau^{w^{\prime}}+1$ let
$p_{\min}=p_{\min}^{\ast}=0$. Therefore $e$ will be either
$\mathbb{P}Qr$ or $\mathbb{P}Q^{\ast}r$. Moreover, notice that every
$\overline{t}$ of type (\ref{eq_bQpq}.a-d) intersects $e$ in a
single vertex. Consequently $\tau^{w}=\tau^{w^{\prime}}+1$. \QED

\begin{Lemma}\label{Rk_3}If $\tau^{w}>0$ and $\overline{t}$ is of type
(\ref{eq_bQpq}.b) or both (\ref{eq_bQpq}.c,d), then:%
\begin{equation}
p=\left\vert \mathbb{P}\right\vert \text{ or }q=\left\vert \mathbb{Q}%
\right\vert \label{eq_pq}%
\end{equation}

\end{Lemma}

\demo If (\ref{eq_pq}) is false then there must be $p_{i}$, $q_{k}$
such that $w\left(  p_{i}\right)  =p_{\min}$, $w\left(
p_{i}^{\ast}\right)  =p_{\min }^{\ast}$, $w\left(  q_{k}\right)
=q_{\min}$, $w\left(  q_{k}^{\ast}\right) =q_{\min}^{\ast}$. Now,
assume $\overline{t}$ is of type (\ref{eq_bQpq}.b), then since
$w\left(  vv^{\ast}\right)  \geq\tau^{w}$ we
have that:%
\[
\tau^{w}=w\left(  p_{i}p_{j}^{\ast}q_{k}q_{l}^{\ast}\right)  =p_{\min}%
+p_{\min}^{\ast}+q_{\min}+q_{\min}^{\ast}=w\left(  p_{i}p_{i}^{\ast}\right)
+w\left(  q_{k}q_{k}^{\ast}\right)  \geq2\tau^{w}%
\]
which is a contradiction.

Now, consider that $\overline{t}$ is of both types
(\ref{eq_bQpq}.c,d), and notice that either $w\left(  r\right)  $ or $w\left(
r^{\ast}\right)  $ are greater than zero (as a consequence of (\ref{eq_bQpq}%
.a)). Assume that $w\left(  r\right)  >0$, then it follows that:%
\[
\tau^{w}=w\left(  p_{i}p_{j}^{\ast}r\right)  =p_{\min}+p_{\min}^{\ast
}+w\left(  r\right)  =w\left(  p_{i}p_{i}^{\ast}\right)  +w\left(  r\right)
\geq\tau^{w}+1
\]
once again, a contradiction. \QED

\begin{Corollary}\label{Cor_3} If $\tau^{w}>0$, $\max\left\{  p_{\min},p_{\min}^{\ast}\right\}
>0$, $\max\left\{  q_{\min},q_{\min}^{\ast}\right\}  >0$ and $\overline{t}$ is of type (\ref{eq_bQpq}.b) or both types (\ref{eq_bQpq}.c,d),
then there is an edge $e$ of $Q_{pq}^{F}$ such that $w^{\prime}\geq0$ and
$\tau^{w}=\tau^{w^{\prime}}+1$.
\end{Corollary}

\demo From Lemma \ref{Rk_3} we know that either $p=\left\vert \mathbb{P}\right\vert \text{ or }q=\left\vert \mathbb{Q}%
\right\vert$. Therefore, consider that $\tau^{w}=w\left(  p_{i}p_{j}^{\ast}%
q_{k}q_{l}^{\ast}\right)  $, then we have that $w\left(  p_{i}p_{j}^{\ast
}q_{k}q_{l}^{\ast}\right)  \leq w\left(  p_{i}p_{j}^{\ast}r\right)  $ and
$w\left(  p_{i}p_{j}^{\ast}q_{k}q_{l}^{\ast}\right)  \leq w\left(  q_{k}%
q_{l}^{\ast}r^{\ast}\right)  $ but this implies that:%
\begin{equation}
w\left(  r\right)  \geq w\left(  q_{k}q_{l}^{\ast}\right)  >0,w\left(
r^{\ast}\right)  \geq w\left(  p_{i}p_{j}^{\ast}\right)  >0 \label{eq_wrr*}%
\end{equation}

On the other hand, if $\tau^{w}=w\left(  p_{i}p_{j}^{\ast}r\right)
=w\left( q_{k}q_{l}^{\ast}r^{\ast}\right)  =\tau^{w}$ we have that
$w\left(  rr^{\ast }\right)  \geq w\left(  p_{i}p_{j}^{\ast}r\right)
,w\left(  q_{k}q_{l}^{\ast }r^{\ast}\right)  $ but this also implies
(\ref{eq_wrr*}). Therefore, we must pick $e$ from (\ref{eq_e_Cor2}).
The selection of $e$ guarantees that $w^{\prime}\geq0$, while the
occurrence of $\mathbb{P}$ or $\mathbb{Q}$ in $e$ guarantees that
$\tau^{w}=\tau^{w^{\prime}}+1$. \QED

\begin{Corollary}\label{Cor_1} If $\tau^{w}>0$, $\max\left\{  p_{\min},p_{\min}^{\ast}\right\}
>0$, $\max\left\{  q_{\min},q_{\min}^{\ast}\right\}  >0$ and $\overline{t}$ is not of type (\ref{eq_bQpq}.b) or both (\ref{eq_bQpq}.c)
or (\ref{eq_bQpq}.d) then there is an $e\in F_{pq}\cup F_{pq}^{\ast}$ such
that $\tau^{w}=\tau^{w^{\prime}}+1$ and $w^{\prime}\geq0$.
\end{Corollary}

\demo Let us assume that $\overline{t} = a_{i}a_{j}^{\ast}c$ is of
type (\ref{eq_bQpq}.c or d). Then $w\left(  cc^{\ast}\right)  \geq
w\left(  a_{i}a_{j}^{\ast}c\right)  $, but this implies that
$w\left( c^{\ast}\right)  \geq w\left(  a_{i}a_{j}^{\ast}\right)
>0$. Therefore $e$
could be picked as:%
\begin{equation}%
\begin{array}
[c]{cc}%
PQc^{\ast} & \text{if }p_{\min},q_{\min}>0\\
PQ^{\ast}c^{\ast} & \text{if }p_{\min},q_{\min}^{\ast}>0\\
P^{\ast}Qc^{\ast} & \text{if }p_{\min}^{\ast},q_{\min}>0\\
P^{\ast}Q^{\ast}c^{\ast} & \text{if }p_{\min}^{\ast},q_{\min}^{\ast}>0
\end{array}
\label{eq_e_Cor1}%
\end{equation}

Moreover, the selection of $e$ guarantees that $w^{\prime}\geq0$ and since
$w\left(  a_{i}a_{j}^{\ast}c\right)  =w^{\prime}\left(  a_{i}a_{j}^{\ast
}c\right)  +1$ and $w\left(  a_{i}a_{j}^{\ast}c\right)  <w\left(  b_{k}%
b_{l}^{\ast}c^{\ast}\right)  $ it is also true that
$w^{\prime}\left( a_{i}a_{j}^{\ast}c\right)  \leq w^{\prime}\left(
b_{k}b_{l}^{\ast}c^{\ast }\right)  +2$. On the other hand, if
$\tau^{w}=w\left(  vv\right)  $ for some $v\in PQr$ then we only
need to pick $e$ in such a way that $w^{\prime}\geq0$ since $w\left(
vv^{\ast}\right)  =w^{\prime}\left(  vv^{\ast }\right)  +1$, and
considering that $c^{\ast}$ is either $r$ or $r^{\ast}$ we can again
pick $e$ from (\ref{eq_e_Cor1}). \QED

As an example, consider $\overline{t} = q_{i}%
q_{j}^{\ast}r^{\ast}$, and $p_{\min},q_{\min}^{\ast}>0$ then from
(\ref{eq_e_Cor1}) we know that $e=PQ^{\ast}r$.

\medskip

By the above results we obtain:

\begin{Theorem}\label{menger2}
The hypergraph $Q_{pq}^{F}$ with $F=F_{pq}^{\ast}$ is mengerian.
\end{Theorem}

\demo By Theorem~\ref{QFpq_Packing} $Q_{pq}^{F}$ has the packing
property. Thus $Q_{pq}^{F}$ is an ideal hypergraph
(Theorem~\ref{lehman}). If $\tau^{\omega} (Q_{pq}^{F}) > 0$ then by
the previous Corollaries \ref{Cor_2}, \ref{Cor_3} and \ref{Cor_1}
there exists $e\in E(Q_{pq}^{F})$ such that $\tau^{\omega} (
Q_{pq}^{F} )=\tau^{\omega -e} ( Q_{pq}^{F} )+1$. Therefore by
Theorem~\ref{menger}, $Q_{pq}^{F}$ is mengerian. \QED

\medskip

Let us denote by $I(Q_{pq}^{F})=\widetilde{I}$ the ideal generated by $F(Q_{pq}^{F})= \{x^{v_1},\ldots,x^{v_k}\}$, where $v_i$ is the $i^{th}$ column of the matrix $A(Q_{pq}^{F})$, then:
%Let $$B(Q_{pq}^{F})=\left(\hspace{-1mm}
%\begin{array}{ccc}
%v_1&\cdots&v_k\\
%1&\cdots &1
%\end{array}\hspace{-1mm}
%\right)
%$$
\begin{Corollary}
$R[\widetilde{I}t]$ is normal and the set covering polyhedron $Q(A(Q_{pq}^{F}))$ is  integral.
\end{Corollary}
\demo It follows by applying Theorem~\ref{noclu1}. \QED

%\begin{Corollary}
%If $B(Q_{pq}^{F})$ has rank $s$, then $\Delta_s(B(Q_{pq}^{F}))=1$.
%\end{Corollary}
%\demo
%By using Theorem~\ref{MM}.
%\QED

\bigskip

\medskip

Finally, we give an algorithm that constructs a list $m$ with $\tau^{w}$ edges from
$Q_{pq}^{F}$ such that $\sum_{e\in m}e\leq w$.

\begin{enumerate}
\item Set $i=0$, $m_{0}=\emptyset$ and $w_{0}=w$

\item \textbf{while} $\tau^{w_{i}}\neq0$

\begin{enumerate}
\item Obtain the values $p_{\min},p_{\min}^{\ast},q_{\min},q_{\min}^{\ast}$
for $w_{i}$

\item \textbf{if} $\max\left\{  p_{\min},p_{\min}^{\ast}\right\}  >0$,
$\max\left\{  q_{\min},q_{\min}^{\ast}\right\}  >0$, \textbf{then}:

\begin{enumerate}
\item \textbf{if} $w_{i}\left(  p_{i}p_{j}^{\ast}r\right)  =w_{i}\left(
q_{k}q_{l}^{\ast}r^{\ast}\right)  =\tau^{w_{i}}$ or $w_{i}\left(  p_{i}%
p_{j}^{\ast}q_{k}q_{l}^{\ast}\right)  =\tau^{w_{i}}$ \textbf{then} pick $e$
from (\ref{eq_e_Cor2})

\item \textbf{else} pick $e$ from (\ref{eq_e_Cor1})
\end{enumerate}

\item \textbf{else} pick $e$ from (\ref{eq_e_Cor2})

\item Set $m_{i+1}=\left(  e,m_{i}\right)  $ and $w_{i+1}=w_{i}-e$

\item $i=i+1$
\end{enumerate}

\item \textbf{return} $m_{i}$
\end{enumerate}

Note that the former algorithm repeats the while statement in line 2, $\tau^{w}$ times. This is so since for every iteration, $e$ is picked in such a way that $\tau^{w_{i}}%
=\tau^{w_{i-1}}+1$. Moreover, if $m$ represents the returned value, we have
that:%
\[
\sum_{e\in m}e_{v}\leq w_{v}\text{ for every }v\in V_{pq}%
\]
since otherwise, the vector $w_{\tau^{w}}$ obtained at the end of the while cycle in line 2
would contain at least one negative entry. Therefore verifying if $Q_{pq}^{F}$ is mengerian. This algorithm can be generalized to other hypergraphs.

\end{document}